\documentclass[12pt]{article}
\usepackage{amssymb,amsmath}
\usepackage{cases}
\usepackage{amsfonts}
\usepackage{color,xcolor,cite}
\usepackage[left=2.0cm,right=2.0cm,top=2.0cm,bottom=2.0cm]{geometry}
\usepackage[colorlinks,citecolor=blue,urlcolor=blue]{hyperref}

\newtheorem{theorem}{Theorem}[section]

\newtheorem{lemma}{Lemma}[section]
\newtheorem{proposition}{Proposition}[section]

\newtheorem{definition}{Definition}[section]
\newtheorem{remark}{Remark}[section]

\newcommand{\bal}{\begin{align}}
\newcommand{\bbal}{\begin{align*}}
\newcommand{\beq}{\begin{equation}}
\newcommand{\eeq}{\end{equation}}
\newcommand{\bca}{\begin{cases}}
\newcommand{\eca}{\end{cases}}

\newcommand{\pa}{\partial}
\newcommand{\fr}{\frac}

\newcommand{\De}{\Delta}

\newcommand{\ep}{\varepsilon}
\newcommand{\dd}{\mathrm{d}}

\newcommand{\R}{\mathbb{R}}
\newcommand{\vv}{\mathbf{v}}
\newcommand{\les}{\lesssim}

\newcommand{\bi}{\Big}

\begin{document}
\title{Sharp ill-posedness for the generalized Camassa-Holm
equation in Besov spaces }

\author{Jinlu Li$^{1}$, Yanghai Yu$^{2,}$\footnote{E-mail: lijinlu@gnnu.edu.cn; yuyanghai214@sina.com(Corresponding author); mathzwp2010@163.com} and Weipeng Zhu$^{3}$\\
\small $^1$ School of Mathematics and Computer Sciences, Gannan Normal University, Ganzhou 341000, China\\
\small $^2$ School of Mathematics and Statistics, Anhui Normal University, Wuhu 241002, China\\
\small $^3$ School of Mathematics and Big Data, Foshan University, Foshan, Guangdong 528000, China}

\date{\today}

\maketitle\noindent{\hrulefill}

{\bf Abstract:} In this paper, we consider the Cauchy problem for the generalized Camassa-Holm
equation that includes the Camassa-Holm as well as the Novikov equation on the line. We present a new and unified method to prove the sharp ill-posedness for the generalized Camassa-Holm
equation in $B^s_{p,\infty}$ with $s>\max\{1+1/p, 3/2\}$ and $1\leq p\leq\infty$ in the sense that the solution map to this equation starting from $u_0$ is
discontinuous at $t = 0$ in the metric of $B^s_{p,\infty}$. Our results cover and improve the previous work given in \cite{Li22}, solving an open problem left in \cite{Li22}.


\vskip0mm\noindent{\hrulefill}

\section{Introduction}\label{sec1}
In this paper, we consider the Cauchy problem for the generalized Camassa-Holm-Novikov (gCHN) equation which was proposed by Anco, Silva and Freire \cite{Anco2015} as follows
\begin{eqnarray}\label{eq1}
        \left\{\begin{array}{ll}
         m_t+u^km_x+(k+1)u^{k-1}u_xm=0,\; &(x,t,k)\in \R\times\R^+\times\mathbb{Z}^+,\\
         m=u-u_{xx},\\
         u(0, x)=u_0(x),\; &x\in \R. \end{array}\right.
        \end{eqnarray}
The gCHN is an evolution equation with (k+1)-order nonlinearities, and can be regarded as a subclass of the generalized Camassa-Holm (g-kbCH) equation considered in \cite{Grayshan2013, Himonas2014}
\begin{equation*}
  m_t+u^km_x+bu^{k-1}u_xm=0,\quad k\in \mathbb{Z}^+,\; b\in\mathbb{R}.
\end{equation*}

When $k=1$, (\ref{eq1}) reduces to the classical Camassa-Holm (CH) equation \cite{Escher1,Escher2,Escher3,Escher4,Escher5}
\begin{eqnarray}\label{eq2}
          u_t-u_{xxt}+3uu_x=2u_xu_{xx}+uu_{xxx},
                  \end{eqnarray}
which was originally derived as a bi-Hamiltonian system by Fokas and Fuchssteiner \cite{Fokas1981} in the context of the KdV model and gained prominence after Camassa-Holm \cite{Camassa1993} independently re-derived
it as an approximation to the Euler equations of hydrodynamics. \eqref{eq2} is completely integrable \cite{Camassa1993,Constantin-P} with a bi-Hamiltonian structure \cite{Constantin-E,Fokas1981} and infinitely many conservation laws \cite{Camassa1993,Fokas1981}. Also, it admits exact peaked
soliton solutions (peakons) of the form $ce^{-|x-ct|}$ with $c>0$, which are orbitally stable \cite{Constantin.Strauss} and models wave breaking (i.e., the solution remains bounded, while its slope becomes unbounded in finite time \cite{Constantin,Escher2,Escher3}).

When $k=2$, (\ref{eq1}) becomes the famous Novikov equation \cite{H-H,HHK2018,Li2,Ni2011,Novikov2009,wu2021,Wu2012,Wu2013,Yan2012}
\begin{eqnarray}\label{eq3}
         u_t-u_{xxt}+4u^2u_x=3uu_xu_{xx}+u^2u_{xxx}.
        \end{eqnarray}
Home-Wang \cite{Home2008} proved that the Novikov equation  with cubic nonlinearity shares similar properties with the CH equation, such as a Lax pair in matrix form, a bi-Hamiltonian structure, infinitely many conserved quantities and peakon solutions given by the formula $u(x, t)=\sqrt{c}e^{-|x-ct|}$.

Setting $\Lambda^{-2}=\big(1-\pa^2_x\big)^{-1}$, then we transform \eqref{eq1} equivalently into the following nonlinear transport type equation
\begin{eqnarray}\label{gch}
        \left\{\begin{array}{ll}
         u_t+u^ku_x=\mathbf{P}(u)+\mathbf{Q}(u),\\
         u(0, x)=u_0(x), \end{array}\right.
        \end{eqnarray}
where
\begin{align}\label{CH1}
\mathbf{P}(u):=-\pa_x\Lambda^{-2}\Big(\frac{2k-1}{2}u^{k-1}u_x^2+u^{k+1}\Big)\quad\text{and}\quad\mathbf{Q}(u):=-\frac{k-1}{2}\Lambda^{-2}\big(u^{k-2}u_x^3\big).
\end{align}
As shown in \cite{Anco2015, Grayshan2013, Himonas2014}, the gCHN equation (\ref{eq1}) admits a local conservation law, possesses single peakons of the form $u(x, t)=c^{1/k}e^{-|x-ct|}$ as well as multi-peakon solutions and exhibits wave breaking phenomena (see \cite{Yan2019} and the references therein).

In recent years, the issue of well-posedness in different spaces for the g-kbCH equation has been a fascinating object of research due to its abundant physical and mathematical properties and a series of achievements have been made in the study of the g-kbCH equation. For general $k$, using a Galerkin-type approximation scheme, Himonas-Holliman \cite{Himonas2014} established the local well-posedness of the g-kbCH equation in the Sobolev space $H^s(\mathbb{R}\;\text{or}\; \mathbb{T})$. Zhao-Li-Yan \cite{Zhao2014} extended the above well-posedness result to the Besov space $B_{p, r}^s(\mathbb{R})$ with  $s>\max\{1+\frac{1}{p}, \frac{3}{2}\}$ and $1\leq p, r\leq \infty$. However, for $r=\infty$, they established the continuity of the data-to-solution map in a weaker topology. Subsequently, Chen-Li-Yan \cite{Chen2015} solved the critical case for $(s, p, r)=(\frac{3}{2}, 2, 1)$. Guo, Liu, Molinet and Yin \cite{Guo-Yin} established the ill-posedness of the Camassa-Holm equation in the critical Sobolev
space $H^{3/2}(\mathbb{R}\;\text{or}\; \mathbb{T})$ and even in the Besov space $B_{p,r}^{1+1/p}(\mathbb{R}\;\text{or}\; \mathbb{T})$ with $p\in[1,\infty],r\in(1,\infty]$ by proving the norm inflation. In our recent paper \cite{Li22}, we proved the solution map to the Camassa-Holm equation starting from $u_0$ is
discontinuous at $t = 0$ in the metric of $B^s_{p,\infty}(\R)$, which implies the ill-posedness for this equation in $B^s_{p,\infty}(\R)$. More precisely, we established
\begin{theorem}[See \cite{Li22}]\label{th0}
Let $s>2+\max\big\{3/2,1+1/p\big\}$ with $1\leq p\leq \infty$. There exits $u_0\in B^s_{p,\infty}(\R)$ and a positive constant $\ep_0$ such that the data-to-solution map $u_0\mapsto \mathbf{S}_{t}(u_0)$ of the Cauchy problem \eqref{gch}-\eqref{CH1} with $k=1$ satisfies
\bbal
\limsup_{t\to0^+}\|\mathbf{S}_{t}(u_0)-u_0\|_{B^s_{p,\infty}}\geq \ep_0.
\end{align*}
\end{theorem}
In addition, for the Novikov equations, we proved that Theorem \ref{th0} holds for only $p=2$ in \cite{Li22}, while for the case $p\neq2$ the difficulty lies mainly in the construction of initial data $u_0$ due to the appearance of $u_0^2\pa_xu_0$. Naturally, the method in \cite{Li22} seems to be invalid when proving the ill-posedness for \eqref{gch}-\eqref{CH1} with $k\geq3$ in $B^s_{p,\infty}(\R)$ ($p\neq2$) since the construction of initial data make the computation of $u_0^k\pa_xu_0$ more difficult.

In this present paper, we shall develop a new and unified method to study the ill-posedness problem for \eqref{gch}-\eqref{CH1} with general $k\in\mathbb{Z}^+$. Our main aim is to prove the solution map to the Cauchy problem \eqref{gch}-\eqref{CH1} starting from $u_0$ is
discontinuous at $t = 0$ in the metric of $B^s_{p,\infty}(\R)$, which implies the ill-posedness for this equation in $B^s_{p,\infty}(\R)$. Furthermore, we expect that the sharp index pair $(s,p)$ satisfies that $s>\max\big\{3/2,1+1/p\big\}$ with $1\leq p\leq \infty$.
Now let us state our main result of this paper.
\begin{theorem}\label{th1}
Let $k\in\mathbb{Z}^+$ be fixed. Assume that
\bal\label{in}
s>\max\Big\{\frac32,1+\frac1p\Big\} \quad\text{with}\quad 1\leq p\leq \infty.
\end{align}
There exits $u_0\in B^s_{p,\infty}(\R)$ and a positive constant $\ep_0$ such that the data-to-solution map $u_0\mapsto \mathbf{S}_{t}(u_0)$ of the Cauchy problem \eqref{gch}-\eqref{CH1}
satisfies
\bbal
\limsup_{t\to0^+}\|\mathbf{S}_{t}(u_0)-u_0\|_{B^s_{p,\infty}}\geq \ep_0.
\end{align*}
\end{theorem}

\begin{remark} As mentioned above, system \eqref{gch}-\eqref{CH1} unifies the Camassa-Holm and Novikov equations. In \cite{Li22}, we  only obtained the ill-posedness for the Novikov equation in $B^s_{2,\infty}$ with $s>\fr72$ due to the technical difficulty. Thus Theorem \ref{th1} covers our recent results on both the Camassa-Holm and Novikov equationsin \cite{Li22}. Also, our method here is new and makes the proof simpler.
\end{remark}
The Cauchy problem for the Degasperis-Procesi equation reads as
\begin{align}\label{dp}
\begin{cases}
\pa_tu+uu_x=-\frac32\pa_x(1-\pa^2_x)^{-1}(u^2), &\quad (t,x)\in \R^+\times\R,\\
u(0,x)=u_0(x), &\quad x\in \R.
\end{cases}
\end{align}
\begin{remark}
We should mention that the $H^1$ norm of solutions to the g-kbCH is conserved if and only if $b=k+1$, which naturally excludes the Degasperis-Procesi equation for the case $k=1$ and $b=3$. However, following the procedure in the proof of Theorem \ref{th1} with suitable modification, we can prove Theorem \ref{th1} holds for the Degasperis-Procesi equation.
\end{remark}

\section{Preliminaries}\label{sec2}
{\bf Notation}\;
The notation $A\les B$ (resp., $A \gtrsim B$) means that there exists a harmless positive constant $c$ such that $A \leq cB$ (resp., $A \geq cB$).
Given a Banach space $X$, we denote its norm by $\|\cdot\|_{X}$. For $I\subset\R$, we denote by $\mathcal{C}(I;X)$ the set of continuous functions on $I$ with values in $X$. Sometimes we will denote $L^p(0,T;X)$ by $L_T^pX$.

Next, we will recall some facts about the Littlewood-Paley decomposition, the nonhomogeneous Besov spaces and some of their useful properties (For more details, see \cite{Bahouri2011}).
Let $\mathcal{B}:=\{\xi\in\mathbb{R}:|\xi|\leq \frac 4 3\}$ and $\mathcal{C}:=\{\xi\in\mathbb{R}:\frac 3 4\leq|\xi|\leq \frac 8 3\}.$
There exist two radial functions $\chi\in C_c^{\infty}(\mathcal{B})$ and $\varphi\in C_c^{\infty}(\mathcal{C})$ both taking values in $[0,1]$ such that
\begin{align*}
&\chi(\xi)+\sum_{j\geq0}\varphi(2^{-j}\xi)=1 \quad \forall \;  \xi\in \R^d.
\end{align*}

\begin{definition}[See \cite{Bahouri2011}]
For every $u\in \mathcal{S'}(\mathbb{R})$, the Littlewood-Paley dyadic blocks ${\Delta}_j$ are defined as follows
\begin{numcases}{\Delta_ju=}
0, & if $j\leq-2$;\nonumber\\
\chi(D)u=\mathcal{F}^{-1}(\chi \mathcal{F}u), & if $j=-1$;\nonumber\\
\varphi(2^{-j}D)u=\mathcal{F}^{-1}\big(\varphi(2^{-j}\cdot)\mathcal{F}u\big), & if $j\geq0$.\nonumber
\end{numcases}
\end{definition}
\begin{definition}[See \cite{Bahouri2011}]
Let $s\in\mathbb{R}$ and $(p,r)\in[1, \infty]^2$. The nonhomogeneous Besov space $B^{s}_{p,r}(\R)$ is defined by
\begin{align*}
B^{s}_{p,r}(\R):=\Big\{f\in \mathcal{S}'(\R):\;\|f\|_{B^{s}_{p,r}(\mathbb{R})}<\infty\Big\},
\end{align*}
where
\begin{numcases}{\|f\|_{B^{s}_{p,r}(\mathbb{R})}=}
\left(\sum_{j\geq-1}2^{sjr}\|\Delta_jf\|^r_{L^p(\mathbb{R})}\right)^{\fr1r}, &if $1\leq r<\infty$,\nonumber\\
\sup_{j\geq-1}2^{sj}\|\Delta_jf\|_{L^p(\mathbb{R})}, &if $r=\infty$.\nonumber
\end{numcases}
\end{definition}
\begin{remark}\label{re3}
It should be emphasized that the fact $B^s_{p,\infty}(\R)\hookrightarrow B^t_{p,\infty}(\R)$ with $s>t$ will be often used implicity.
\end{remark}
Finally, we give some important properties which will be also often used throughout the paper.
\begin{lemma}[See \cite{Bahouri2011}]\label{le1}
Let $(p,r)\in[1, \infty]^2$ and $s>\max\big\{1+\frac1p,\frac32\big\}$. Then we have
\bbal
&\|uv\|_{B^{s-2}_{p,r}(\R)}\leq C\|u\|_{B^{s-2}_{p,r}(\R)}\|v\|_{B^{s-1}_{p,r}(\R)}.
\end{align*}
\end{lemma}

\begin{lemma}[See \cite{Bahouri2011}]\label{le2}
For $(p,r)\in[1, \infty]^2$, $B^{s-1}_{p,r}(\R)$ with $s>1+\frac{1}{p}$ is an algebra. Moreover, for any $u,v \in B^{s-1}_{p,r}(\R)$ with $s>1+\frac{1}{p}$, we have
\bbal
&\|uv\|_{B^{s-1}_{p,r}(\R)}\leq C\|u\|_{B^{s-1}_{p,r}(\R)}\|v\|_{B^{s-1}_{p,r}(\R)}.
\end{align*}
\end{lemma}
\begin{remark}\label{re5} Let $(p,r)\in[1, \infty]^2$ and $s>\max\big\{1+\frac1p,\frac32\big\}$, using Lemmas \ref{le1}-\ref{le2}, we have
\begin{itemize}
  \item for the terms $\mathbf{P}(u)$ and $\mathbf{P}(v)$, there holds
\begin{align}
\|\mathbf{P}(u)-\mathbf{P}(v)\|_{B_{p, r}^{s-1}} & \lesssim\|u-v\|_{B_{p, r}^{s-1}}\big(\|u\|_{B_{p, r}^{s-1}}+\|v\|_{B_{p, r}^{s-1}}\big)^{k-1}\big(\|u\|_{B_{p, r}^{s}}+\|v\|_{B_{p, r}^{s}}\big).
\end{align}
  \item for the terms $\mathbf{Q}(u)$ and $\mathbf{Q}(v)$ (notice that $\mathbf{Q}(u)=0$ for $k=1$), there holds
\begin{align}
        \|\mathbf{Q}(u)-\mathbf{Q}(v)\|_{B_{p, r}^{s-1}}\lesssim
         \|u-v\|_{B_{p, r}^{s-1}}\big(\|u\|_{B_{p, r}^{s-1}}+\|v\|_{B_{p, r}^{s-1}}\big)^{k-2}\big(\|u\|_{B_{p, r}^{s}}+\|v\|_{B_{p, r}^{s}}\big)^2.
\end{align}
\end{itemize}
\end{remark}
\begin{lemma}[See \cite{Bahouri2011}]\label{le3}
For $1\leq p\leq \infty$ and $s>0$. There exists
a constant $C$, depending continuously on $p$ and $s$, swe have
\bbal
\bi\|2^{j s}\left\|[\Delta_j,v]\pa_xf\right\|_{L^{p}}\bi\|_{\ell^{\infty}} \leq C\big(\|\pa_x v\|_{L^{\infty}}\|f\|_{B_{p, \infty}^{s}}+\|\pa_x f\|_{L^{\infty}}\|\pa_xv\|_{B_{p, \infty}^{s-1}}\big),
\end{align*}
where we denote the standard commutator $[\Delta_j,v]\pa_xf=\Delta_j(v\pa_xf)-v\Delta_j\pa_xf$.
\end{lemma}

\section{Proof of Theorem \ref{th1}}\label{sec3}
\subsection{Construction of Initial Data}\label{sec3.1}
We need to introduce smooth, radial cut-off functions to localize the frequency region. Precisely,
let $\widehat{\phi}\in \mathcal{C}^\infty_0(\mathbb{R})$ be an even, real-valued and non-negative function on $\R$ and satisfy
\begin{numcases}{\widehat{\phi}(\xi)=}
1,&if $|\xi|\leq \frac{1}{4}$,\nonumber\\
0,&if $|\xi|\geq \frac{1}{2}$.\nonumber
\end{numcases}

\begin{lemma}\label{le4}
Define the function $f_n(x)$ by
$$f_n(x)=\phi(x)\cos \Big(\frac{17}{12}2^{n}x\Big)\quad\text{with}\quad n\gg1.$$
Then we have
\begin{numcases}{\Delta_j(f_n)=}
f_n, &if $j=n$,\nonumber\\
0, &if $j\neq n$.\nonumber
\end{numcases}
\end{lemma}
{\bf Proof.}\; See \cite{Li1}.

\begin{lemma}\label{le5}
Define the initial data $u_0(x)$ as
\bbal
u_0(x):=\sum\limits^{\infty}_{n=0}2^{-ns} \phi(x)\cos \bi(\frac{17}{12}2^{n}x\bi).
\end{align*}
Then for any $s>\max\{\frac32,1+\frac1p\}$ and $k\in\mathbb{Z}^+$, we have for some $n$ large enough
\bbal
&\|u_0\|_{B^{s}_{p,\infty}}\leq C,\\
&\|u^k_0\pa_x\De_{n}u_0\|_{L^p}\geq c2^{n(1-s)},
\end{align*}
where $C$ and $c$ are some positive constants.
\end{lemma}
{\bf Proof.}\; By the definition of Besov space and the support of $\varphi(2^{-j}\cdot)$, we have
\bbal
\|u_0\|_{B^{s}_{p,\infty}}
&\leq C.
\end{align*}
Using Lemma \ref{le4} yields
\bbal
\De_{n}u_0(x)&=2^{-ns} \phi(x)\cos \bi(\frac{17}{12}2^{n}x\bi),
\end{align*}
equivalently,
\bbal
\pa_x\De_{n}u_0&=2^{-ns} \phi'(x)\cos \bi(\frac{17}{12}2^{n}x\bi)-\frac{17}{12}2^{n}2^{-ns} \phi(x)\sin \bi(\frac{17}{12}2^{n}x\bi).
\end{align*}
Thus, we have
\bbal
u^k_0\pa_x\De_{n}u_0&=2^{-ns} u^k_0(x)\phi'(x)\cos \bi(\frac{17}{12}2^{n}x\bi)-\frac{17}{12}2^{n}2^{-ns} u^k_0(x)\phi(x)\sin \bi(\frac{17}{12}2^{n}x\bi).
\end{align*}
Since $u^k_0(x)$ is a real-valued and continuous function on $\R$, then there exists some $\delta>0$ such that
\bal\label{yh}
&|u^k_0(x)|\geq \fr{1}{2}|u^k_0(0)|=\fr{1}{2}\Big(\phi(0)\sum\limits^{\infty}_{n=0}2^{-ns}\bi)^k=\frac{2^{sk}\phi^k(0)}{2(2^s-1)^k}\quad\text{ for any }  x\in B_{\delta}(0).\end{align}
Thus we have from \eqref{yh}
\bbal
\|u^k_0\pa_x\De_{n}u_0\|_{L^p}
&\geq C2^{n}2^{-ns} \bi\|\phi(x)\sin \bi(\frac{17}{12}2^{n}x\bi)\bi\|_{L^p(B_{\delta}(0))}-C2^{-ns}\bi\| \phi'(x)\phi^{k}(x)\cos \bi(\frac{17}{12}2^{n}x\bi)\bi\|_{L^p}\\
&\geq (c2^{n}-C)2^{-ns} .
\end{align*}
We choose $n$ large enough such that $C<\frac{c}{2}2^{n}$ and then finish the proof of Lemma \ref{le5}.
\subsection{Error Estimates}\label{sec3.2}
\begin{proposition}\label{pro3.1}
Assume that $\|u_0\|_{B^s_{p,\infty}}\lesssim 1$. Let $u\in L^\infty_TB^{s}_{p,\infty}$ be the solution of the Cauchy problem \eqref{gch}, then under the assumptions of Theorem \ref{th1}, we have
\begin{align}\label{w}
&\|\mathbf{S}_{t}\left(u_{0}\right)-u_0\|_{B^{s-1}_{p,\infty}}\lesssim t.
\end{align}
Furthermore, there holds
\begin{align}\label{w1}
\left\|\mathbf{w}\right\|_{B_{p, \infty}^{s}} \lesssim t^{2},
\end{align}
here and in what follows we denote
$$\mathbf{w}:=\mathbf{S}_{t}(u_0)-u_0-t\mathbf{v}_0 \quad \text{with}\quad\mathbf{v}_0:=\mathbf{P}(u_0)+\mathbf{Q}(u_0)-u_0^k\pa_x u_0.$$
\end{proposition}
{\bf Proof.}\; For simplicity, we denote $u(t):=\mathbf{S}_t(u_0)$ here and in what follows. Due to the fact $B^s_{p,\infty}\hookrightarrow\rm Lip$, we know that there exists a positive time $T=T(\|u_0\|_{B^s_{p,\infty}})$ such that
\bbal
\|u(t)\|_{L^\infty_TB^s_{p,\infty}}\leq C\|u_0\|_{B^s_{p,\infty}}\leq C.
\end{align*}
Using the Mean Value Theorem and Remark \ref{re5} with $v=0$, we obtain from \eqref{gch} that
\bal\label{s}
\|u(t)-u_0\|_{B^{s-1}_{p,\infty}}
&\leq \int^t_0\|\pa_\tau u\|_{B^{s-1}_{p,\infty}} \dd\tau
\nonumber\\&\leq \int^t_0\big(\|\mathbf{P}(u)\|_{B^{s-1}_{p,\infty}}+\int^t_0\|\mathbf{Q}(u)\|_{B^{s-1}_{p,\infty}}\big) \dd\tau+ \int^t_0\|u^k\pa_xu\|_{B^{s-1}_{p,\infty}} \dd\tau
\nonumber\\&\les t\|u\|^{k+1}_{L_t^\infty B^{s}_{p,\infty}}
\nonumber\\&\les t\|u_0\|^{k+1}_{B^{s}_{p,\infty}}\nonumber\\&\les t.
\end{align}
By the Mean Value Theorem and \eqref{gch}, then we obtain that
\begin{align}\label{p1}
\|\mathbf{w}\|_{B^{s-2}_{p,\infty}}
&\leq \int^t_0\|\partial_\tau u-\mathbf{v}_0\|_{B^{s-2}_{p,\infty}} \dd\tau \nonumber\\
&\les \int^t_0\big(\|\mathbf{P}(u)-\mathbf{P}(u_0)\|_{B^{s-2}_{p,\infty}}+\|\mathbf{Q}(u)-\mathbf{Q}(u_0)\|_{B^{s-2}_{p,\infty}}\big) \dd\tau \nonumber\\&\quad+\int^t_0\|u^k\partial_xu-u^k_0\partial_xu_0\|_{B^{s-2}_{p,\infty}}\dd\tau.
\end{align}
Using Remark \ref{re5} again yields
\begin{align}\label{p2}
\int^t_0\big(\|\mathbf{P}(u)-\mathbf{P}(u_0)\|_{B^{s-2}_{p,\infty}}+\|\mathbf{Q}(u)-\mathbf{Q}(u_0)\|_{B^{s-2}_{p,\infty}}\big) \dd\tau \les t^2,
\end{align}
Notice that the simple fact $$a^{k+1}-b^{k+1}=(k+1)(a-b)\xi^k,\quad \xi\; \text{is between}\; a \;\text{and}\; b,$$
we obtain from Lemma \ref{le2} and \eqref{s} that
\begin{align}\label{p3}
\|u^k\partial_xu(\tau)-u^k_0\partial_xu_0\|_{B^{s-2}_{p,\infty}}&\les \|u^{k+1}(\tau)-u_0^{k+1}\|_{B^{s-1}_{p,\infty}}\nonumber\\
&\les\|u(\tau)-u_0\|_{B^{s-1}_{p,\infty}} \|u_0\|^k_{B^{s-1}_{p,\infty}}
\nonumber\\
&\les \tau.
\end{align}
Inserting \eqref{p2} and \eqref{p3} into \eqref{p1} yields
\begin{align*}
\|\mathbf{w}\|_{B^{s-2}_{p,\infty}}\les t^2.
\end{align*}
Thus, we finish the proof of Proposition \ref{pro3.1}.

Now we present the proof of Theorem \ref{th1}.\\
{\bf Proof of Theorem \ref{th1}.}\quad Notice that $$\mathbf{S}_{t}(u_0)-u_0=t\mathbf{v}_0+\mathbf{w}(t,u_0) \quad \text{and}\quad\mathbf{v}_0=\mathbf{P}(u_0)+\mathbf{Q}(u_0)-u_0^k\pa_x u_0,$$
by the triangle inequality and Propositions \ref{pro3.1}, we deduce that
\bal\label{M}
\|\mathbf{S}_{t}(u_0)-u_0\|_{B^{s}_{p,\infty}}
&\geq2^{{ns}}
\big\|\De_{n}\big(\mathbf{S}_{t}(u_0)-u_0\big)\big\|_{L^p}\nonumber\\
&=2^{{ns}}\big\|\De_{n}\big(t\vv_0+\mathbf{w}(t,u_0)\big)\big\|_{L^p}
\nonumber\\&\geq t2^{{ns}}\|\De_{n}\big(\mathbf{v}_0\big)\|_{L^p}
-2^{{2n}}2^{{n(s-2)}}
\big\|\De_{n}\big(\mathbf{w}(t,u_0)\big)\big\|_{L^p}\nonumber\\
&\geq t2^{{n}s}\|\De_{n}\big(u_0^k\pa_xu_0\big)\|_{L^p}-
t2^{{n}s}\|\De_{n}\big(\mathbf{P}(u_0)+\mathbf{Q}(u_0)\big)\|_{L^p}\nonumber\\
&\quad-C2^{2{n}}\|\mathbf{w}(t,u_0)\|_{B^{s-2}_{p,\infty}}
\nonumber\\&\geq t2^{{n}s}\|u_0^k\pa_x\De_{n}u_0\|_{L^p}-t2^{{n}s}\|[\De_{n},u_0^k]\pa_xu_0\|_{L^p}\nonumber\\
&\quad-
t\|\mathbf{P}(u_0)+\mathbf{Q}(u_0)\|_{B^{s}_{p,\infty}}-C2^{2{n}}t^2\nonumber\\&\geq t2^{{n}s}\|u_0^k\pa_x\De_{n}u_0\|_{L^p}-Ct\big\|2^{{n}s}\|[\De_{n},u_0^k]\pa_xu_0\|_{L^p}\big\|_{\ell^\infty}\nonumber\\
&\quad-t\|\mathbf{P}(u_0)+\mathbf{Q}(u_0)\|_{B^{s}_{p,\infty}}-C2^{2{n}}t^2.
\end{align}
By Lemmas \ref{le2}-\ref{le3}, one has
\bbal
\|\mathbf{P}(u_0)+\mathbf{Q}(u_0)\|_{B^{s}_{p,\infty}}\les \|u_0^{k-1}(\pa_xu_0)^2+u_0^{k+1}\|_{B^{s-1}_{p,\infty}}+
\|u_0^{k-2}(\pa_xu_0)^3\|_{B^{s-1}_{p,\infty}}\les 1
\end{align*}
and
\bbal
\big\|2^{{n}s}\|[\De_{n},u_0^k]\pa_xu_0\|_{L^p}\big\|_{\ell^\infty}\les \|\pa_x(u_0^k)\|_{L^\infty}\|u_0\|_{B^{s}_{p,\infty}}+
\|\pa_xu_0\|_{L^\infty}\|\pa_x(u_0^k)\|_{B^{s-1}_{p,\infty}}\les 1.
\end{align*}
Gathering all the above estimates and Lemma \ref{le5} together with \eqref{M}, we obtain
\bbal
\|\mathbf{S}_{t}(u_0)-u_0\|_{B^{s}_{p,\infty}}\geq ct2^{{n}}-Ct-C2^{2{n}}t^2.
\end{align*}
Taking large $n$ such that $c2^{{n}}\geq 2C$, we have
\bbal
\|\mathbf{S}_{t}(u_0)-u_0\|_{B^{s}_{p,\infty}}\geq \frac{c}{2}t2^{{n}}-C2^{2{n}}t^2.
\end{align*}
Thus, picking $t2^{n}\approx\ep$ with small $\ep$, we have
\bbal
\|\mathbf{S}_{t}(u_0)-u_0\|_{B^{s}_{p,\infty}}\geq \frac{c}{2}\ep-C\ep^2\geq c_1\ep.
\end{align*}
This completes the proof of Theorem \ref{th1}.

\vspace*{1em}
\noindent\textbf{Acknowledgements.} J. Li is supported by the National Natural Science Foundation of China (11801090 and 12161004) and Postdoctoral Science Foundation of China (2020T130129 and 2020M672565). Y. Yu is supported by the National Natural Science Foundation of China (12101011) and Natural Science Foundation of Anhui Province (1908085QA05). W. Zhu is partially supported by the National Natural Science Foundation of China (11901092) and Natural Science Foundation of Guangdong Province (2017A030310634).

\vspace*{1em}
\noindent\textbf{Conflict of interest}
The authors declare that they have no conflict of interest.

\end{document}